\documentclass[review]{elsarticle}
\usepackage{lineno}

\usepackage[T1]{fontenc} 
\usepackage{microtype} 
\usepackage{bm,subcaption}
\usepackage{amsmath,amssymb,amsthm}
\usepackage[english]{babel} 

\usepackage[right=3cm,top=3cm,bottom=3cm,left=3cm]{geometry} 
\usepackage{pgfplots,tikz} 
\pgfplotsset{compat=1.13}
\usepackage{algorithm}
\usepackage[noend]{algpseudocode}
\usepackage{booktabs,multirow,multicol}
\usepackage[export]{adjustbox}
\usepackage{tikz}
\usepackage{mathabx,color,url}
\usepackage[colorlinks=true]{hyperref}

\usepackage{lineno}
\modulolinenumbers[1]
\usepackage{wrapfig}
\usepgfplotslibrary{statistics}
\usetikzlibrary{pgfplots.statistics} 
\usepgfplotslibrary{groupplots}

\allowdisplaybreaks

\makeatletter
\def\munderbar#1{\underline{\sbox\tw@{$#1$}\dp\tw@\z@\box\tw@}}
\makeatother

\pgfplotsset{compat=1.12}

\usepackage{enumitem} 
\newcommand{\rev}{\textcolor{black}}

\begin{document}
\hypersetup{
	breaklinks=true,
	urlcolor=blue, 
	linkcolor=blue, 
	citecolor=red, 
	pdftitle={},
	pdfauthor={},
}

\begin{frontmatter}

\title{\Large \textbf{Piecewise Polyhedral Formulations for a Multilinear Term}} 

\author[lanl]{Kaarthik Sundar}
\author[lanl]{Harsha Nagarajan\corref{corr_author}}
\cortext[corr_author]{Corresponding author}
\ead{harsha@lanl.gov}
\author[wisconsin]{Jeff Linderoth}
\author[clemson]{Site Wang}
\author[lanl]{Russell Bent}

\address[lanl]{Los Alamos National Laboratory, Los Alamos, NM, USA}
\address[wisconsin]{Department of Industrial and Systems Engineering, Department of Computer Sciences, University of Wisconsin-Madison, Madison, WI, USA}
\address[clemson]{Department of Industrial Engineering, Clemson University, SC, USA}

\date{} 

\begin{abstract}
In this paper, we present a mixed-integer linear programming (MILP) formulation of a piecewise, polyhedral relaxation (PPR) of a multilinear term using its convex-hull representation. Based on the PPR's solution, we also present a MILP formulation whose solutions are feasible for nonconvex, multilinear equations. We then present computational results showing the effectiveness of proposed formulations on standard benchmark nonlinear programs (NLPs) with multilinear terms and compare with a traditional formulation that is built using recursive bilinear groupings of multilinear terms. 
\end{abstract}



\end{frontmatter}

\section{Introduction} \label{sec:intro}
In the global optimization literature, many solution techniques used to solve nonconvex Mixed-Integer Nonlinear Programs (MINLPs) or Nonlinear Programs (NLPs) to optimality rely on convex relaxations \cite{Tawarmalani2005} and piecewise convex relaxations of nonconvex functions \cite{Hasan2010,Deassis2017,Castro2015,Nagarajan2016,Nagarajan2017} to obtain tight bounds on the optimal objective values of MINLPs/NLPs. The three primary topics which this letter focuses are as follows. \textit{First}, on developing piecewise, polyhedral relaxations (PPR) of a multilinear term using the convex hull representation and \textit{second}, on comparing this relaxation with traditional, recursive bilinear piecewise relaxations. Our work generalizes formulations that approximate functions with piecewise, linear functions \cite{Sridhar2013,Vielma2010,Lee2001,Padberg2000} to formulations that relax functions with piecewise, convex polyhedra. \textit{Third}, since recovering feasible solutions for MINLPs/NLPs is of great interest for quantifying the quality of the solutions of PPR, this letter focuses on developing a MILP formulation which requires solutions to satisfy multilinear terms, thereby producing a locally optimal solution that is feasible for all nonconvex, multilinear equations.

Throughout the rest of the letter, boldface is used to denote vectors and we formally define a multilinear term as $\phi(\bm x): [\bm \ell, \bm u] \rightarrow \mathbb R$, where 
\begin{flalign}
\phi(\bm x) = \prod_{i \in \mathcal I} x_i. \label{eq:multilinear_def}
\end{flalign}

\noindent Here, $\mathcal I$ is an index set for the set of variables, $\bm x$, and $[\bm \ell, \bm u] = \{\bm x \in \mathbb R^{|\mathcal I|} : \bm \ell \leqslant \bm x \leqslant \bm u\}$. For every $i \in \mathcal I$, we associate a spatial disjunction with variable $x_i$ defined by discretization points $\mathcal S_i = \{s_{i,1}, s_{i,2}, \dots, s_{i,n_i}\}$ and associated intervals $\mathcal P_i = \{ \delta_i^1, \delta_i^2, \cdots, \delta_i^{n_i-1} \}$ where the intervals form a partition of the domain of $x_i$ given by $[\ell_i, u_i]$, i.e. $\delta_i^k = [s_{i,k}, s_{i, k+1}]$ and $\ell_i = s_{i,1} < s_{i,2} < \cdots < s_{i,n_i} = u_i$. Given this notation, each $x_i$, $i \in \mathcal I$ is constrained to satisfy
\begin{flalign}
x_i \in [s_{i,1}, s_{i,2}]  ~\lor ~ [s_{i,2}, s_{i,3}] ~ \lor ~ \cdots ~ \lor ~ [s_{i,{n_i-1}}, s_{i,n_i}]. \label{eq:spatial-disjunction}
\end{flalign}
Given $\mathcal P_i$ and $\mathcal S_i$ for every $i\in \mathcal I$, we let $(n_i-1)$ and $n_i$ be the cardinality of the respective sets. For any $i \in \mathcal I$ and $k \in \{1, \cdots, n_i-1\}$, the interval $\delta_i^k$ is  defined to be \textit{\text{active}} when the value of $x_i$ lies in the interval $\delta_i^k$. Eq. \eqref{eq:spatial-disjunction} then enforces one interval to be active per variable. Finally, the full set of discretization points and the partition sets defined by these points and intervals in the space $[\bm \ell, \bm u]$ are given by the sets $\mathcal S = \bigtimes_{i \in \mathcal I}\mathcal S_i$ and $\mathcal P = \bigtimes_{i \in \mathcal I}\mathcal P_i$, respectively.

Given the multilinear term $\phi(\bm x)$ and the associated sets $\mathcal S$ and $\mathcal P$, this work presents piecewise, polyhedral relaxations (PPR) for the graph of $\phi(\bm x)$, given by
\begin{flalign}
X = \{ (\bm x, w) \in [\bm \ell, \bm u] \times \mathbb R : w = \phi(\bm x) \}. \label{eq:graph_phi}
\end{flalign}

An important special case occurs when $\phi(\bm x)$ is bilinear with a single partition on each variable, \textit{i.e.}, $|\mathcal I| = 2$ and $|\mathcal P|=1$. Here, the McCormick relaxation, given by the following equations  
\begin{subequations}
\begin{flalign}
w \geqslant u_2 x_1 + u_1 x_2 - u_1 u_2, &\quad w \geqslant \ell_2 x_1 + \ell_1 x_2 - \ell_1 \ell_2,  \label{eq:under} \\
w \leqslant u_2 x_1 + \ell_1 x_2 - \ell_1 u_2, &\quad w \leqslant \ell_2 x_1 + u_1 x_2 - u_1 \ell_2, \label{eq:over}
\end{flalign}
\label{eq:mccormick}
\end{subequations}
\noindent
defines the convex hull \cite{McCormick1976,Al-Khayyal1983} of the set $X$, i.e., $\operatorname{conv}(X)$. For general multilinear terms \textit{i.e.}, $|\mathcal I| \geqslant 3$ with $|\mathcal P|=1$, McCormick described a procedure that recursively applies the McCormick relaxation to products of variables. We refer to this procedure as the ``recursive McCormick relaxation''.
On general multilinear terms with asymmetric bounds on the variables, despite not capturing the convex hull \cite{Luedtke2012}, the recursive McCormick relaxation for a multilinear term is the basis for relaxations used in many global optimization solvers such as BARON \cite{Sahinidis1996}, SCIP \cite{vigerske2018scip} and Couenne \cite{Belotti2009}. 
In contrast, a relaxation based on a vertex representation, for $|\mathcal I| \geqslant 2$ and $|\mathcal P| = 1$ does characterize the convex hull of $X$ \cite{Floudas2013,Rikun1997} and is given by
\begin{flalign}
\operatorname{conv}(X) = \operatornamewithlimits{Proj}_{\bm x,w}~~\Bigl\{ (\bm x, w, \bm \lambda) \in [\bm \ell, \bm u] \times \mathbb R \times \Delta_{2^{|\mathcal I|}} : \bm x = \sum_{s=1}^{2^{|\mathcal I|}} \lambda_s \hat {\bm x}_s, \,w = \sum_{s=1}^{2^{|\mathcal I|}} \lambda_s \phi(\hat{\bm x}_s) \Bigr\}, 
\label{eq:conv}
\end{flalign}
where $\hat {\bm x}_1, \hat {\bm x}_2, \dots, \hat {\bm x}_{2^{|\mathcal I|}}$ is the collection of all the vertices of the hyper-rectangle $[\bm \ell, \bm u]$ given by $\mathcal S$
and $\Delta_{2^{|\mathcal I|}}$ is the $2^{|\mathcal I|}$-dimensional 0-1 simplex. 

\rev{Note that for the case where $|\mathcal I| = 2$ and $|\mathcal P| = 1$ (bilinear term without any partitions), the set of feasible points are identical for formulations in Eq. \eqref{eq:conv} and Eq. \eqref{eq:mccormick}, and thus are equivalent. When $|\mathcal I| \geqslant 2$ and $|\mathcal P| \geqslant 2$, the resulting MILP relaxation is expressed as a disjunctive union of $|\mathcal P|$ polytopes where each polytope is obtained by applying Eq. \eqref{eq:conv} to the domain defined by the corresponding partition set}. The convex hull of this disjunctive union is completely characterized by utilizing the theory of disjunctive programming \cite{Balas1979,Jeroslow1984} and the introduction of a sufficient number of auxiliary variables; in particular, by utilizing the formulations in \cite{Jeroslow1984}. The formulation in \cite{Jeroslow1984} is an extended formulation that generates constraints for each partition separately and then aggregating them. In contrast, the approach presented in this letter works directly with the combinatorial structure underlying the shared extreme points and presents non-extended PPR for a multilinear term. \rev{Furthermore, logarithmically-sized MILP relaxations for disjunctive constraints in the special case of a bilinear term with spatial disjunction on one variable is addressed in \cite{Misner2012,Misener2011}, and it's generalization to multilinear terms in \cite{huchette2019combinatorial}, albeit without much computational studies. To the best of our knowledge, the work in \cite{Misner2012,Misener2011,huchette2019combinatorial} is the state-of-the-art for PPR of a multilinear term.}

The contributions of this letter include (a) a special-order-set (SOS)-2-based PPR that uses a vertex representation for a single multilinear term and a spatial disjunction on each variable, and (b) given the PPR's solution, a piecewise formulation that uses the vertex representation of the convex hull polytope to recover the feasible solution for the nonconvex multilinear equation. This PPR and the feasible solutions obtained are compared with a relaxation and feasible solutions recovered based on recursively relaxing bilinear groupings of the multilinear term akin to the recursive McCormick relaxation. We refer to this relaxation as the recursive piecewise polyhedral relaxation (R-PPR). To the best of our knowledge, this is the first work that presents a systematic comparison between an SOS-2-based PPR and a recursive piecewise relaxation, along with the feasible solution recovery for a general multilinear term with spatial disjunctions on every variable. 

\section{Piecewise Polyhedral Relaxation for a Multilinear Term} \label{sec:ppr}
In this section, we present the PPR for $X$ given $\mathcal I$ (the index set of variables), $\mathcal P$ (the partition set), and $\mathcal S$ (the set of discretization points). For notational simplicity, for any $i \in \mathcal{I}$ and $k \in \{1, \dots, n_i - 1\}$, we use $\munderbar{\delta}_i^k = s_{i,k}$ and $\widebar{\delta}_i^k = s_{i,k+1}$ to denote the lower and upper limits for this interval, respectively. We associate a non-negative multiplier variable, $\lambda_{s}$, with each $\hat{\bm x}_s \in \mathcal S$. We also define $y_{i,k}$ as a binary indicator for the $k$\textsuperscript{th} interval of variable $x_i$, $i\in \mathcal I$. This variable takes a value of $1$ when interval $\delta_i^k$ is active for variable $x_i$ and takes a value of $0$ otherwise. We let $\bm \lambda$ and $\bm y$ denote the vector of continuous multiplier variables and indicator variables, respectively. Next, for each $i \in \mathcal I$, we define a set-valued function, $\mu_i(r) = \{\hat{\bm x}_s \in \mathcal S: \bm e_i^T \hat{\bm x}_s = r\}$ where $\bm e_i$ is a unit vector whose $i$\textsuperscript{th} coordinate is equal to $1$ and is $0$ everywhere else. This function defines the subset of points in $\mathcal S$ whose $i$\textsuperscript{th} component is equal to $r$.
Finally, we let $\lambda(\hat{\mathcal S}) = \sum_{\hat{\bm x}_s \in \hat{\mathcal S}} \lambda_s$ for any 
$\hat{\mathcal S} \subseteq \mathcal S$.

Given these notations, the PPR for a multilinear term is given by 
\begin{equation}
    \left \langle \prod_{i \in \mathcal I} x_i \right \rangle^{\text{PPR}} = 
    \left \{ 
    (\bm x, w, \bm \lambda, \bm y) \in [\bm \ell, \bm u] \times \mathbb R \times \Delta_{|\mathcal S|} \times \{0,1\}^{\sum_{i \in \cal{I}} (n_i-1)} \,\, \left|  \,\,
    (\bm x, w, \bm \lambda, \bm y) \text{ satisfies Eq. \eqref{eq:sos2} }
    \right. \right \} \label{eq:ppr}
\end{equation}
where, Eq. \eqref{eq:sos2} is given by
\begin{subequations}\label{ppr}
\begin{flalign}
&\bm x = \sum_{\hat{\bm x}_s \in \mathcal S} \lambda_s \hat{\bm x}_s, \, w = \sum_{\hat{\bm x}_s \in \mathcal S} \lambda_s \phi(\hat{\bm x}_s), &\label{eq:xzl} \\
& \sum_{k=1}^{n_i-1} y_{i,k} = 1, \quad \forall i \in \mathcal I, & \label{eq:y} \\
& \lambda(\mu_i(\munderbar{\delta}_i^1)) \leqslant y_{i,1}, \quad \forall i \in \mathcal I, & \label{eq:sos2-1} \\
& \lambda(\mu_i(\munderbar{\delta}_i^k)) \leqslant y_{i,{k-1}}+y_{i,k}, \quad \forall i \in \mathcal I,\,\, k \in 2,\dots,n_i-2, & \label{eq:sos2-2} \\
& \lambda(\mu_i(\munderbar{\delta}_i^{n_i-1})) \leqslant y_{i,n_i-1}, \quad \forall i \in \mathcal I, & \label{eq:sos2-3} \\
&\lambda_s \geqslant 0, \quad \forall \hat{\bm x}_s \in \mathcal S, \text{ and, } &\label{eq:lambdann} \\
& y_{i,k} \in \{0,1\}, \quad \forall i \in \mathcal I, \,\, k \in 1,2,\dots,n_i-1. & \label{eq:ybin}
\end{flalign}
\label{eq:sos2}
\end{subequations}

\noindent
Constraints in Eq. \eqref{ppr} use the binary variables of each variable partition to activate and deactivate the non-negative multiplier variables $\lambda_s$, $s\in \mathcal S$. In particular, Eq. \eqref{eq:y} ensures that only one partition per variable is active. The 
constraints in Eq. \eqref{eq:sos2-1} -- \eqref{eq:sos2-3} enforce adjacency conditions on the $\lambda_s$ variables akin to SOS-2 constraints \cite{Beale1976a, Beale1976}. \rev{In the context of mathematical models for formulating general disjunctive union of polytopes, this model is also referred to as the ``Convex-Combination (CC)" model or the ``lambda" formulation in the literature \cite{Vielma2010}. The 
constraints in Eq. \eqref{eq:sos2-1} -- \eqref{eq:sos2-3} ensure that the relaxation is the convex combination of the vertices of the active partition in $\mathcal P$.} The traditional SOS-2 is an ordered set of non-negative variables, of which at most two can be non-zero and they must be consecutive in their ordering. In this case, instead of an ordered set of non-negative variables, we have an ordered set that is a sum of non-negative variables \emph{i.e.}, for any $i\in \mathcal I$, $\lambda(\mu_i(\munderbar{\delta}_i^k)))$, with $k=\{1,\dots, n_i-1\}$, forms an ordered set of non-negative variable sums where SOS-2 constraints are imposed. 

It is not difficult to establish that the PPR in Eq. \eqref{eq:ppr} has the property of being locally sharp, i.e., the projection of the LP relaxation of Eq. \eqref{eq:ppr} to $(\bm x, w)$-space is exactly the $\operatorname{conv}(X)$. This is a property of the SOS-2 formulation for piecewise, linear approximations of functions \cite{Huchette2017}. 
In the literature, there are other vertex-based formulations which use this variable space and other formulations which use a logarithmic number of binary variables \cite{Huchette2017,Padberg2000,Misener2012global}, however, we
focus on providing a detailed comparison study between the SOS-2-based PPR in \eqref{eq:sos2} and its recursive counterparts (in Sec. \ref{sec:r-ppr}). 

\paragraph{PPR Example}
The PPR for a bilinear term with $|\mathcal P|=9$ ($3$ partitions on each variable), shown in Fig. \ref{fig:setup}, is given in Eq. \eqref{eq:f1_ex}.
\begin{figure}
\centering
\begin{tikzpicture}[scale=0.8]
\draw [<->] (5,0) node [below] {$x_1$} -- (0,0) -- (0,5) node [left] {$x_2$};
\draw [thick] (1,1) node [left] {$\hat{\bm x}_1$} -- (2,1) node [below] {$\hat{\bm x}_2$} -- (3,1) node [below] {$\hat{\bm x}_3$} -- (4,1) node [right] {$\hat{\bm x}_4$};
\draw [thick] (1,2) node [left] {$\hat{\bm x}_5$} -- (2,2) node [below left] {$\hat{\bm x}_6$} -- (3,2) node [below right] {$\hat{\bm x}_7$} -- (4,2) node [right] {$\hat{\bm x}_8$};
\draw [thick] (1,3) node [left] {$\hat{\bm x}_9$} -- (2,3) node [above left] {$\hat{\bm x}_{10}$} -- (3,3) node [above right] {$\hat{\bm x}_{11}$} -- (4,3) node [right] {$\hat{\bm x}_{12}$};
\draw [thick] (1,4) node [left] {$\hat{\bm x}_{13}$} -- (2,4) node [above] {$\hat{\bm x}_{14}$} -- (3,4) node [above] {$\hat{\bm x}_{15}$} -- (4,4) node [right] {$\hat{\bm x}_{16}$};
\draw [thick] (1,1) -- (1,4);
\draw [thick] (2,1) -- (2,4);
\draw [thick] (3,1) -- (3,4);
\draw [thick] (4,1) -- (4,4);
\node [below] at (1,0) {$s_{1,1}$}; \draw[fill] (1,0) circle [radius=0.025];
\node [below] at (2,0) {$s_{1,2}$}; \draw[fill] (2,0) circle [radius=0.025];
\node [below] at (3,0) {$s_{1,3}$}; \draw[fill] (3,0) circle [radius=0.025];
\node [below] at (4,0) {$s_{1,4}$}; \draw[fill] (4,0) circle [radius=0.025];
\node [left] at (0,1) {$s_{2,1}$}; \draw[fill] (0,1) circle [radius=0.025];
\node [left] at (0,2) {$s_{2,2}$}; \draw[fill] (0,2) circle [radius=0.025];
\node [left] at (0,3) {$s_{2,3}$}; \draw[fill] (0,3) circle [radius=0.025];
\node [left] at (0,4) {$s_{2,4}$}; \draw[fill] (0,4) circle [radius=0.025];
\end{tikzpicture}
\caption{$\phi(\bm x) = x_1 x_2$ with $|P_1| = |P_2| = 3$}
\label{fig:setup}
\end{figure}
\vspace{-0.2cm}
\begin{subequations}
\begin{flalign}
& \bm x = \sum_{s=1}^{16} \lambda_s \hat{\bm x}_s, \, w = \sum_{s=1}^{16} \lambda_s \phi(\hat{\bm x}_s), \, \sum_{s=1}^{16} \lambda_s = 1, & \\
& y_{1,1} + y_{1,2} + y_{1,3} = 1, \,\, y_{2,1} + y_{2,2} + y_{2,3} = 1, & \\
& \lambda_1 + \lambda_5 + \lambda_9 + \lambda_{13} \leqslant y_{1,1}, \quad \quad \quad \quad \  \lambda_1 + \lambda_2 + \lambda_3 + \lambda_{4} \leqslant y_{2,1}, &  \\
& \lambda_{2} + \lambda_{6} + \lambda_{10} + \lambda_{14} \leqslant y_{1,1} + y_{1,2}, \quad
\lambda_5 + \lambda_6 + \lambda_7 + \lambda_{8} \leqslant y_{2,1} + y_{2,2}, &  \\
& \lambda_{3} + \lambda_{7} + \lambda_{11} + \lambda_{15} \leqslant y_{1,2} + y_{1,3}, \quad \lambda_{9} + \lambda_{10} + \lambda_{11} + \lambda_{12} \leqslant y_{2,2} + y_{2,3}, &  \\
& \lambda_{4} + \lambda_{8} + \lambda_{12} + \lambda_{16} \leqslant y_{1,3}, \quad \quad \quad \quad
\lambda_{13} + \lambda_{14} + \lambda_{15} + \lambda_{16} \leqslant y_{2,3}, &  \\
& \lambda_s \geqslant 0 \quad \forall s = 1,\dots,16 \quad \text{ and}, \quad y_{i,k} \in \{0,1\} \quad \forall i=1,2, \, k=1,2,3.  
\end{flalign}
\label{eq:f1_ex}
\end{subequations}
\vspace{-0.8cm}
\section{Recursive Piecewise Polyhedral Relaxation for a Multilinear Term with $|\mathcal I| \geqslant 3$} 
\label{sec:r-ppr}
For a general multilinear term with $|\mathcal I| \geqslant 3$ and no partitions i.e., $|\mathcal P| = 1$, the recursive McCormick relaxation successively uses the McCormick relaxation in Eq. \eqref{eq:mccormick} on bilinear groupings of the terms. We extend this approach to build a R-PPR by recursively grouping bilinear terms and applying the PPR, given in Eq. \eqref{eq:ppr}, to each bilinear term. For a particular grouping of bilinear terms with $|\mathcal I| \geqslant 3$, the R-PPR is succinctly expressed as 

\begin{equation}
    \left \langle \prod_{i \in \mathcal I} x_i \right \rangle^{\text{R-PPR}} = 
    \left \langle \left \langle \left \langle \left \langle x_1 \cdot x_2 \right \rangle^{\text{PPR}} \dots \right \rangle^{\text{PPR}}  \cdot x_{|\mathcal I|-1} \right \rangle^{\text{PPR}} \cdot x_{|\mathcal I|}\right \rangle^{\text{PPR}} 
    \label{eq:r-ppr}
\end{equation}

Notice that a different R-PPR can be created by choosing a different order of recursive grouping of variables, leading to different relaxation qualities \cite{belotti2013composition,lee2018algorithmic}. In computational studies, we present the performance of different groupings of bilinear terms in addition to the grouping presented in \eqref{eq:r-ppr}. Also, while applying R-PPR, only the original space of variables in the multilinear term are partitioned, and not the lifted/auxiliary variables resulting from recursive relaxations.

\section{MILP-based Piecewise Formulation for Recovering Feasible Solutions} 
\label{sec:feasible} 
\rev{Thus far, the focus of the paper has been on developing piecewise polyhedral relaxations for multilinear terms. The optimal solution to this relaxation provides an active partition for each variable. We let $\mathcal P^*$ and $\mathcal S^*$ denote the optimal active partition and the extreme points of the active partition chosen by the the PPR of the multilinear term, respectively. Let $[\bm \ell^*, \bm u^*]$ denote the variable bound vector defined by the active partition $\mathcal P^*$. Applying the convex hull formulation given in Eq. \eqref{eq:conv} for $\bm x \in [\bm \ell^*, \bm u^*]$ yields the same solution produced by the PPR. Additionally, enforcing the solution to lie on the one-dimensional faces of the polytope in Eq. \eqref{eq:conv} for variable bounds $[\bm \ell^*, \bm u^*]$ yields a feasible solution to the multilinear term, which is straight-forward to observe. This provides a simple MILP-based technique to compute a feasible solution of the original nonlinear problem of interest. Hence, in this section, we present an alternate MILP that recovers a feasible solution by requiring the solution to lie on an edge of the polytope which represents the relaxation of the multilinear term's active partition. The formulation presented below is easily generalized to the case of requiring a solution to lie on the edges of the convex hull relaxation for any partition, but we omit this for the sake of simplicity of exposition. To enforce these additional constraints, we introduce some additional notation. Let $\mathcal V$ and $\mathcal E$ denote the vertices and the edges (one-dimensional faces) of the polytope defined by Eq. \eqref{eq:conv} for $x \in [\bm \ell^*, \bm u^*]$. We know that $|\mathcal{V}| = 2^{|I|}$ and $|\mathcal{E}| = |I|\cdot 2^{|I|-1}$. Then, for every $v \in \mathcal V$, we define a set function $\gamma(v) = \{e \in \mathcal E: \text{ $e$ is incident on $v$} \}$. Additionally, we introduce a binary variable $z_e$ for $e \in \mathcal E$ that takes a value $1$ if the solution lies on the edge $e$. Given these notations, the piecewise formulation that enforces the solution to lie on one of the edges of the polytope defined by Eq. \eqref{eq:conv} for $x \in [\bm \ell^*, \bm u^*]$ is defined by}
\begin{equation}
    \left \{ 
    (\bm x, w, \bm \lambda, \bm z) \in [\bm \ell^*, \bm u^*] \times \mathbb R \times \Delta_{|\mathcal{V}|} \times \{0,1\}^{|\mathcal{E}|} \,\, \left|  \,\,
    (\bm x, w, \bm \lambda, \bm z) \text{ satisfies Eq. \eqref{eq:feas3} }
    \right. \right \} \label{eq:feas3-short}
\end{equation}

\noindent
where, Eq. \eqref{eq:feas3} is given by
\begin{subequations}
\label{eq:feas1}
\begin{flalign}
&\bm x = \sum_{\hat{\bm x}_s \in \mathcal S^*} \lambda_s \hat{\bm x}_s, \quad w = \sum_{\hat{\bm x}_s \in \mathcal S^*} \lambda_s \phi(\hat{\bm x}_s), &\label{eq:xzl-fs} \\
& \sum_{e \in \mathcal{E}} z_e = 1, \quad \lambda_v \leqslant \sum_{e \in \gamma(v)} z_e \ \forall v \in \mathcal{V}, &\label{eq:main-fs} \\
&\lambda_s \geqslant 0, \quad \forall \hat{\bm x}_s \in \mathcal S^*, \text{ and, } &\label{eq:lambdann-fs} \\
& z_e \in \{0,1\} \quad \forall e \in \mathcal E, & \label{eq:zbin-fs}
\end{flalign}
\label{eq:feas3}
\end{subequations} 
\noindent
where Eqs. \eqref{eq:xzl-fs} and \eqref{eq:main-fs} ensure that any solution to the formulation \eqref{eq:feas3} is a convex combination of the vertices that constitute a one-dimensional face of the polytope in Eq. \eqref{eq:conv} and that only one face is active. \rev{Note that the above described method is also applicable if the original problem is an MINLP, as recovering feasible solutions to MINLPs is by itself an arduous task using off-the-shelf solvers}. In the next section, we present extensive computational results that compare PPR with R-PPR, and the quality of recovered feasible solutions.

\section{Computational Results}
In this section, we present three sets of computational results. The first two sets of results use the NLP instance shown in Eq. \eqref{eq:instance1} and the third set of results uses 60 nonlinear optimization problems that contain only multilinear terms taken from \cite{Bao2015}. 
\begin{subequations}
\begin{flalign}
\mathcal F = \max  & \quad (x_1 x_2 x_3 x_4 + x_3 x_4 x_5 x_6 + x_5 x_6 x_7 x_8) & \\
\text{subject to:}  & \quad  100 x_1 - x_2 - x_3 + 833 x_4 + 95 x_5 + x_6 - x_7 + 100x_8 \leqslant 50000, \label{eq:instance1_1}\\
& \quad  100 \leqslant x_1 \leqslant 500, \  1000 \leqslant x_2,x_3 \leqslant 2000, \ 10 \leqslant x_i \leqslant 100 \quad \forall i \in \{4,\dots,8\}. \label{eq:instance1_2}
\end{flalign}
\label{eq:instance1}
\end{subequations}
The first set of results shows the strength of the various relaxations (PPR and R-PPRs) and the effectiveness of feasible solution recovery performed on the active partition obtained from the respective relaxations. The second set of results is aimed at comparing the different ways in which feasible solution recovery can be performed by using the template formulation in Sec. \ref{sec:feasible} in terms of computation time and quality of the feasible solutions. The final set of results compares the strength of the proposed PPRs on instances used in the literature for which the optimal solutions are known a-priori. For both the first and the third set of computational experiments, we partition the variables of the multilinear terms with \textit{uniformly} located discretization points.

All models were implemented in Julia using JuMP v0.19.5 \cite{dunning2017jump} and the experimental results used Gurobi v8.0 \cite{gurobi} with Gurobi's presolver and heuristics turned off. \rev{All nonconvex models were solved using BARON v19.12.7 with CPLEX 12.10 and Ipopt 3.12.8 as BARON's sub-solvers}. The computational experiments were performed on an Intel Xeon CPU L5420 @2.50GHz with 64 GB RAM.

\subsection{Effectiveness of the formulations on the NLP in \eqref{eq:instance1}} \label{subsec:set-1}
Table \ref{tab:UB_LB_gaps} compares the objective value of the PPR with three R-PPRs (denoted by upper bound, $UB$), each with a different variable grouping, on the NLP in Eq. \eqref{eq:instance1}. Note that the choice of grouping changes the quality of the relaxation. The active partition of each variable in the relaxed solution is used to compute a feasible solution (whose objective is denoted by lower bound, $LB$) through formulation Eq. \eqref{eq:feas3}. In this table, $UB$ and $LB$ gaps are given by $\frac{UB-OPT}{UB}\cdot100$ and $\frac{OPT-LB}{LB}\cdot100$, respectively; $OPT$ corresponds to the global optimum value of \eqref{eq:instance1} which is equal to $3.2642E10$. It is clear from these results that, in both the relaxation and the feasible solution recovery formulations, PPR outperforms all the R-PPRs. For example, with only 4 partitions, PPR finds a relaxed solution ($UB$) that is better than any R-PPR (up to 12 partitions). It is also clear that the quality of the feasible solutions ($LB$) recovered based on the active partition of the PPR formulation is consistently better than those of R-PPR's for every chosen number of partitions. Note that, the $LB$ gaps \textit{do not} necessarily decrease monotonically as the number of partitions increase; the new partitions for PPR/R-PPR formulations yield different active partitions, thus leading to different local optimal solutions.
In these experiments, the run-times of both PPR and R-PPR formulations were less than a minute and hence we do not report them explicitly. As is the case with any piecewise formulation, the computational performance of PPR degrades as the number of partitions increases on larger NLPs, however the strength of PPR's relaxation remains superior compared to any other recursive variation.
\begin{table}[H]
    \centering
    \rev{\begin{tabular}{c|ccccccc}
        \toprule 
        \textbf{\# of partitions} & \multicolumn{1}{c}{\multirow{2}{*}{\textbf{2}}} & \multicolumn{1}{c}{\multirow{2}{*}{\textbf{4}}} & \multicolumn{1}{c}{\multirow{2}{*}{\textbf{6}}} & \multicolumn{1}{c}{\multirow{2}{*}{\textbf{8}}} & \multicolumn{1}{c}{\multirow{2}{*}{\textbf{10}}} & 
        \multicolumn{1}{c}{\multirow{2}{*}{\textbf{12}}}\\
        \textbf{per variable} &  &  & &  & &  \\
        \midrule 
        \textbf{\% gaps} & $LB$, $UB$ & $LB$, $UB$ & $LB$, $UB$ & $LB$, $UB$ & $LB$, $UB$ & $LB$, $UB$ \\
        $x_a x_b x_c x_d$ & 2.33, 23.99 & 0.15, 3.20 & 1.11, 2.98 & 0.15, 0.83 & 0.00, 0.69 & 0.05, 0.43\\
        $\langle x_a \langle x_b x_c \rangle \rangle x_d$ & 2.33, 65.47 & 36.74, 25.37 & 2.33, 21.73 & 0.15, 14.72 & 0.00, 12.98 & 2.33, 10.34 \\
        $\langle \langle x_a x_b \rangle x_c  \rangle x_d$ & 2.33, 65.47 & 36.74, 25.37 & 2.33, 21.73 & 0.15, 14.72 & 0.00, 12.98 & 2.33, 10.34 \\ 
        $ x_a \langle x_b \langle x_c  x_d \rangle \rangle$ & 2.33, 47.37 & 5.73, 25.27 & 1.11, 16.78  & 0.15, 12.29 & 0.00, 9.59 & 1.11, 8.19\\
        \bottomrule
    \end{tabular}}
    \caption{Relative optimality gaps of upper and lower bounds in percent for the PPR and the R-PPRs for varying number of partitions per variable. The $\langle \cdot \rangle$ represents the grouping of the different variables into bilinear terms in the R-PPRs.}
    \label{tab:UB_LB_gaps}
\end{table}
\noindent
\textit{BARON's performance on the NLP in \eqref{eq:instance1}}: In order to compare the relative optimality gaps obtained by the PPR formulation in Table \ref{tab:UB_LB_gaps}, we also solved instance \eqref{eq:instance1} using BARON solver. \rev{The best optimality gap BARON reports after a run time of one hour is \textit{20.36\%} with 153,555 open spatial branch-and-bound nodes, where the best upper bound reported by BARON is $4.0986 \times 10^{10}$ and the best lower bound by Ipopt within BARON is $3.2642\times 10^{10}$}. In contrast, the PPR reports \textit{3.34\%} with only \textit{four} partitions per variable (see table \ref{tab:UB_LB_gaps}). The best gap PPR reports is \textit{0.47\%} (12 partitions) which is close to global optimum, thus reinforcing the value of the convex hull formulations proposed in this letter. 

\subsection{Feasible solution recovery for the NLP in \eqref{eq:instance1}} \label{subsec:set-2}
Here, we discuss the quality of the feasible solutions obtained by using different variants of the MILP formulation presented in Sec. \ref{sec:feasible} for recovering solutions to the NLP described in \eqref{eq:instance1}. To this end, we introduce auxiliary variables $t_1,t_2$ and $t_3$, such that the objective function of \eqref{eq:instance1} will be to maximize $(t_1+t_2+t_3)$ subject to constraints $t_1 = x_1x_2x_3x_4, t_2 = x_3x_4x_5x_6, t_3 = x_5x_6x_7x_8$, \eqref{eq:instance1_1} and \eqref{eq:instance1_2}. 
The MILP formulation, as presented in Sec. \ref{sec:feasible}, when applied to each multilinear term in Eq. \eqref{eq:instance1}, results in a feasible solution to the NLP. This formulation computes the best solution that lies on the one-dimensional face of the active partition of each partitioned variable as dictated by the solution of the relaxation. From here on, we refer to this formulation as $\mathcal F^{a}$. Instead, the formulation can be extended to compute the best solution that lies on the one-dimensional face of any combination of the variable partition (the formulation is a trivial extension of the formulation in Eq. \eqref{eq:feas3}, and is not presented). We denote this formulation as $\mathcal F^f_1$. 
Though formulation $\mathcal F^f_1$ may sometimes lead to better quality  feasible solutions, the number of one-dimensional faces grows exponentially with the number of partitions on each variable and thus can result in a MILP with a large number of binary variables. Note that while applying $\mathcal F^f_1$ on various recursive grouping of bilinear terms, we do not partition the auxiliary lifted variables. The third formulation, $\mathcal F_2^f$, is applicable only when using recursive groupings of bilinear terms. Here, we enforce partitions on the auxiliary lifted variables introduced while performing the recursion. We also remark that comparing the solutions of $\mathcal F^a$ and $\mathcal F^f_1$ with $\mathcal F_2^f$ is not a fair comparison. $\mathcal F_2^f$ naturally provides better quality solutions since more variables are partitioned. Nevertheless, we compare all the formulations in Table \ref{tab:fs-recovery}. In this table, all percent gaps are computed relative to the global optimum value ($OPT$), that is, $\frac{OPT-LB}{LB}\cdot100$.

\begin{table}[]
    \centering
    \begin{tabular}{c|cc|ccc|ccc|ccc}
        \toprule 
        \textbf{\# of partitions} & \multicolumn{2}{|c|}{$x_a x_b x_c x_d$} & \multicolumn{3}{|c|}{$\langle x_a \langle x_b x_c \rangle \rangle x_d$} & \multicolumn{3}{|c|}{$\langle \langle x_a x_b \rangle x_c  \rangle x_d$} & \multicolumn{3}{|c}{$ x_a \langle x_b \langle x_c  x_d \rangle \rangle$} \\
        \cmidrule(lr){2-12} 
        \textbf{per variable} & $\mathcal F^a$ & $\mathcal F^f_1$ & $\mathcal F^a$ & $\mathcal F^f_1$ & $\mathcal F^f_2$ & $\mathcal F^a$ & $\mathcal F^f_1$ & $\mathcal F^f_2$ & $\mathcal F^a$ & $\mathcal F^f_1$ & $\mathcal F^f_2$ \\
        \midrule
        \textbf{2} & 2.33 & 2.33 & 86.25 & 2.33	& 0.17 & 86.25 & 2.33 & 0.17 & 2.33 & 2.33 & 0.01 \\ 
        \textbf{4} & 0.15 & 0.15 & 0.15 & 36.74 & 0.05 & 0.15 & 36.74 & 0.01 & 2.33 & 5.73 & 2E-4 \\
        \textbf{6} & 1.11 & 1.11 & 4.22 & 2.33 & 0.03 & 4.22 & 2.33 & 0.03 & 1.11 & 1.11 & 4E-3 \\
        \textbf{8} & 0.15 & 0.15 & 0.15 & 0.15 & 5E-3 & 0.15 & 0.15 & 0.05 & 0.17 & 0.15 & 2E-5 \\
        \textbf{10} & 0.0 & 0.0 & 1.12 & 0.0 & 5E-4 & 1.12 & 0.0 & 5E-3 & 0.0 & 0.0 & 4E-5 \\
        \textbf{12} & 0.05 & 0.05 & 0.15 & 2.33 & 6E-5 & 0.15 & 2.33 & 6E-5 & 0.05 & 1.11 & 2E-5 \\
        \bottomrule
    \end{tabular}
    \caption{Quality (percent gaps) of feasible solutions obtained using different variants of formulation \eqref{eq:feas3} on the NLP in Eq. \eqref{eq:instance1}.}
    \label{tab:fs-recovery}
\end{table}
As expected and consistent with our previous observations, even in the feasible solution recovery process, the convex hull representation of the multilinear terms outperforms the recursive approaches when comparing formulations $\mathcal F^a$ and $\mathcal F^f_1$. As for formulation $\mathcal F^f_2$, it outperforms even the convex hull formulation as it partitions the auxiliary lifted variables. As a side note, it is not difficult to prove that for a fixed partition count, $\mathcal F^f_2$ on any recursive bilinear grouping will always produce the best feasible solution among all the formulations i.e., $\mathcal F^f_1$ and $\mathcal F^a$. 

\subsection{Strength of the relaxations on instances from \cite{Bao2015}} \label{subsec:set-3}
\begin{figure*}[htp]
    \centering
    \begin{subfigure}[t]{\textwidth}
        \centering
        \includegraphics[scale=0.63]{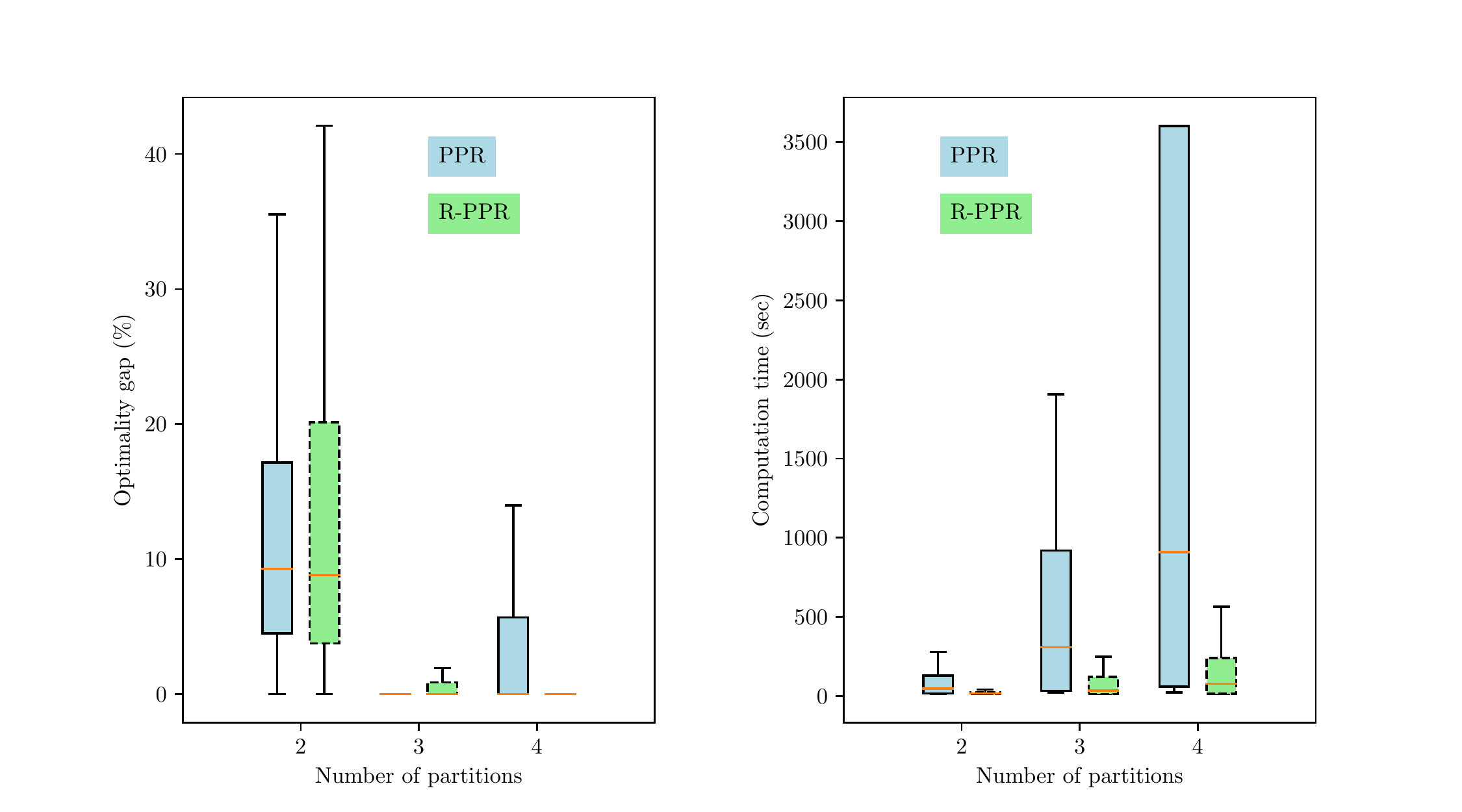}
        \caption{Percent gaps and run times for instances with variable bounds set to $[0,1]$. The optimality gaps for the PPRs increase from $3$ partitions to $4$ partitions because of the computation time limit of $1$ hour imposed on every run. }
    \end{subfigure}%
    \\
    \begin{subfigure}[t]{\textwidth}
        \centering
        \includegraphics[scale=0.63]{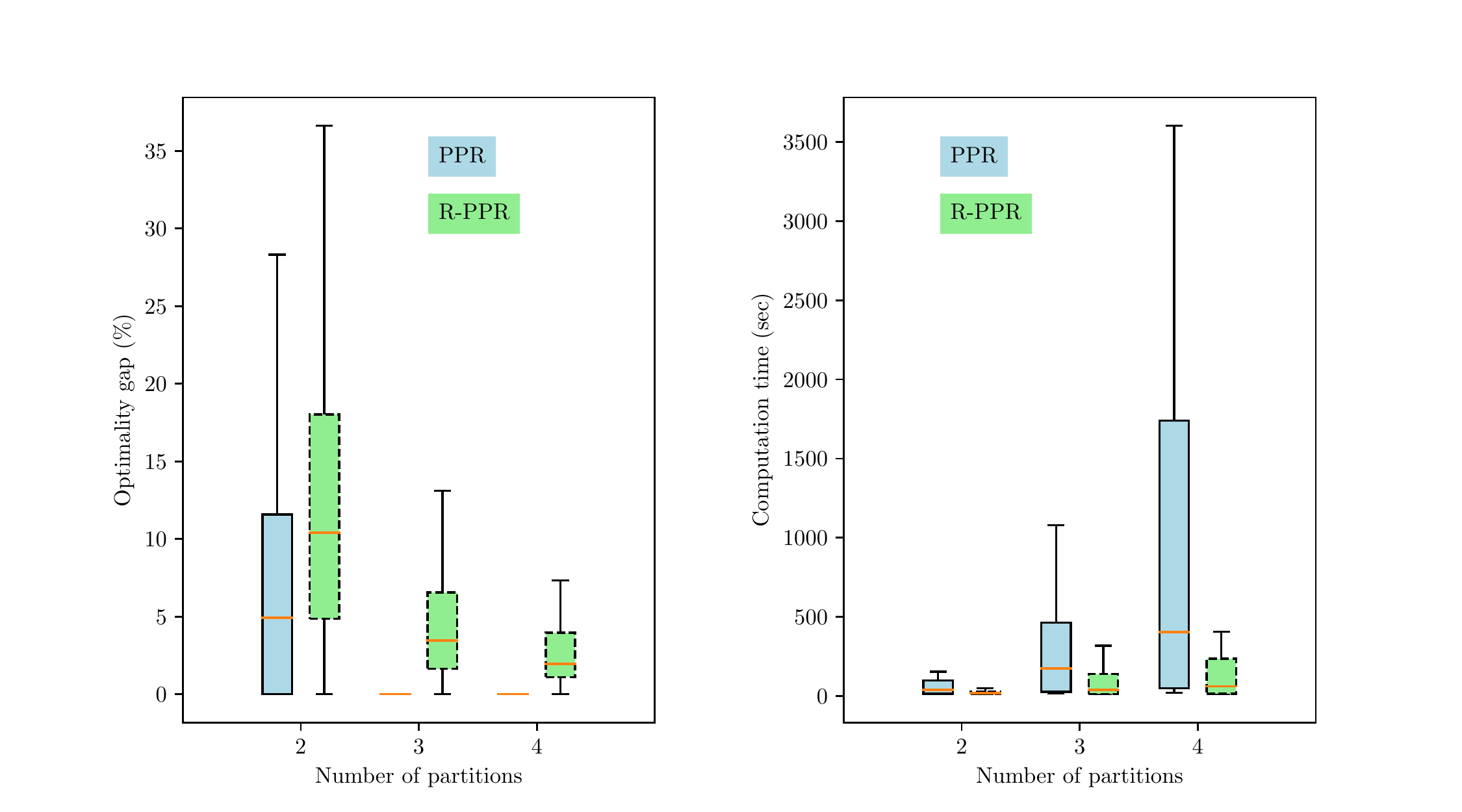}
        \caption{Percent gaps and run times for instances with variable bounds randomly chosen as described earlier. In general, the R-PRR has better computation times because it has fewer multiplier variables ($\lambda$). Nevertheless, the PPR closes the optimality gap with fewer partitions despite taking a little more computation time.}
    \end{subfigure}
    \caption{The optimality gap in \% refers to the relative gap between the objective value of the relaxation of the problem and the global optimum to the problem instance. The computation time is the time taken to solve the PPR and the R-PPR to optimality.}
    \label{fig:box_plot}
\end{figure*}

Our final set of results are based on a collection of 60 multilinear problems taken from \cite{Bao2015}, which are available at: \href{https://minlp.com/downloads/testlibs/barlibs/mult3.zip}{https://minlp.com/downloads/testlibs/barlibs/mult3.zip}. The bounds on all the variables for every problem instance in \cite{Bao2015} are $[0,1]$. To demonstrate the strength of the relaxations discussed in this paper, we created a modified set of instances. For each variable, we uniformly distributed the lower and upper bounds in the interval $[0.1,0.2]$ and $[0.9,1]$, respectively. This modification is motivated by the following known result in \cite{Luedtke2012}--the recursive McCormick relaxation for a multilinear term is the convex hull when all variable bounds are symmetric about the origin or are $[0,1]$. Though the bounds were modified, the global optimum values for all the instances remained the same as those in \cite{Bao2015} when solved using BARON. For each problem instance we relaxed the problems using PPR and a R-PPR with a lexicographic ordering. A time limit of one hour was imposed on every run of the instance. The results of all the runs are summarized using the box-plot shown in Fig. \ref{fig:box_plot}. 

In Fig. \ref{fig:box_plot}, the optimality gap for the PPR in most instances is $0\%$ when the number of partitions per variable is \textit{three}. This is not the case for R-PPR. However, as expected, Fig. \ref{fig:box_plot}(a) suggests that when the bound on every variable in the problem is $[0,1]$, R-PPR is generally the better alternative. 
This is not surprising since similar recursive relaxations are known to capture the convex hull of a multilinear term with symmetric variable bounds \cite{Luedtke2012}. Once the variable bounds are randomized and made asymmetric, the PPR formulation outperforms the recursive relaxation (R-PPR) in both optimality gap and computation time (Fig. \ref{fig:box_plot}(b)). Though the (bottom right) figure suggests that the computation times of PPR are larger than that of R-PPR, it is worth noting that the number of partitions and the run times necessary for the PPR to achieve 0\% optimality gap on all 60 instances is much superior when compared to the R-PPR. 

\section{Concluding Remarks}
This letter presents piecewise polyhedral relaxation formulations for a multilinear term with spatial disjunctions on every variable. An SOS-2-based formulation that generalizes similar formulations for piecewise linear approximations of functions was developed to formulate the convex hull of a single multilinear term. Extensive computational experiments demonstrated that our proposed formulations have significant computational advantages over the well-known recursive grouping of bilinear terms. A MILP-based piecewise formulation was developed to recover high-quality feasible solutions for problems with multilinear terms. Given exponential growth in the number of variables in the PPR formulation, future work will focus on developing formulations with subsets of variables for those rare cases when the number of variables in a multilinear term is large and further add valid inequalities, by generalizing the ideas proposed in \cite{Bao2015}. 

\section*{Acknowledgements}
The work was funded by the Center for Nonlinear Studies (CNLS)
at LANL and LANL's Directed Research and Development (LDRD) projects, ``20170201ER: POD: A Polyhedral Outer-approximation, Dynamic-discretization optimization  solver" and ``20190590ECR: Discrete Optimization Algorithms for Provably Optimal Quantum Circuit Design". This work was carried out under the U.S. DOE Contract No. DE-AC52-06NA25396.

\section{References}
\bibliographystyle{elsarticle-num}
\bibliography{pch.bib}

\end{document}